\documentclass[letterpaper]{jpconf}
\usepackage{amsfonts}
\usepackage{graphicx}
\usepackage{url}

\newcommand{\RR}{{\mathbb R}}
\newcommand{\CC}{{\mathbb C}}
\newcommand{\HH}{{\mathbb H}}
\newcommand{\OO}{{\mathbb O}}
\newcommand{\OP}{{\mathbb {OP}}}
\newcommand{\ZZ}{{\mathbb Z}}

\newcommand{\II}{{\bf I}}
\newcommand{\MM}{{\bf M}}
\newcommand{\PP}{{\bf P}}
\newcommand{\XX}{{\bf X}}

\newcommand{\AAA}{{\cal A}}
\newcommand{\BBB}{{\cal B}}
\newcommand{\III}{{\cal I}}
\newcommand{\MMM}{{\cal M}}
\newcommand{\PPP}{{\cal P}}
\newcommand{\VVV}{{\cal V}}
\newcommand{\XXX}{{\cal X}}

\newcommand{\shat}{\hat{s}}
\newcommand{\that}{\hat{t}}

\newcommand{\eplus}{e_{\scriptscriptstyle\uparrow}}
\newcommand{\eminus}{e_{\scriptscriptstyle\downarrow}}

\newcommand{\phd}{^{\phantom\dagger}}

\renewcommand{\bar}{\overline}
\renewcommand{\tilde}{\widetilde}

\begin{document}

\title{Octonions, \boldmath$E_6$, and Particle Physics}

\author{Corinne A Manogue$^1$ and Tevian Dray$^2$}

\address{$^1$ Department of Physics, Oregon State University,
Corvallis, OR  97331, USA.}
\address{$^2$ Department of Mathematics, Oregon State University,
Corvallis, OR  97331, USA.}

\ead{corinne@physics.oregonstate.edu, tevian@math.oregonstate.edu}

\begin{abstract}
In 1934, Jordan et al.\ gave a necessary algebraic condition, the Jordan
identity, for a sensible theory of quantum mechanics.  All but one of the
algebras that satisfy this condition can be described by Hermitian matrices
over the complexes or quaternions.  The remaining, exceptional Jordan algebra
can be described by $3\times3$ Hermitian matrices over the octonions.

We first review properties of the octonions and the exceptional Jordan
algebra, including our previous work on the octonionic Jordan eigenvalue
problem.  We then examine a particular real, noncompact form of the Lie group
$E_6$, which preserves determinants in the exceptional Jordan algebra.

Finally, we describe a possible symmetry-breaking scenario within $E_6$:
first choose one of the octonionic directions to be special, then choose
one of the $2\times2$ submatrices inside the $3\times3$ matrices to be
special.  Making only these two choices, we are able to describe many
properties of leptons in a natural way.  We further speculate on the ways
in which quarks might be similarly encoded.
\end{abstract}

\section{Introduction}

A personal note from Corinne: During the academic year 1986/87, Tevian and I
were living in York, newly married and young postdocs.  Tevian was working in
the mathematics department there, doing research in general relativity, and I
was working in Durham, with David Fairlie, just beginning my research into the
octonionic structures reported here.  Imagine my pleasure, when I found out
that York had its own resident expert on the octonions!  I returned to York
the following summer, to work with Tony on an attempt to describe the
superstring using octonions~\cite{Sudbery}.  I will be forever grateful to
him, not only for the generous way in which he shared his vast knowledge and
experience in this field, but also for the friendship, respect, collegiality,
and mentorship, which he also generously shared.

\section{Exceptional Quantum Mechanics}

In the Dirac formulation of quantum mechanics, a quantum mechanical state is
represented by a \textit{complex} vector $v$, often written as $|v\rangle$,
which is usually normalized such that $v^\dagger v=1$.  In the Jordan
formulation~\cite{Jordan,JNW,GT,Other}, the same state is instead represented
by the Hermitian matrix $vv^\dagger$, also written as $|v\rangle\langle{v}|$,
which squares to itself and has trace $1$.  The matrix $vv^\dagger$ is thus
the projection operator for the state $v$, which can also be viewed as a pure
state in the density matrix formulation of quantum mechanics.  Note that the
usual phase freedom in $v$ is no longer present in $vv^\dagger$, which is
uniquely determined by the state (and the normalization condition).

A fundamental object in the Dirac formalism is the probability amplitude
$v^\dagger w$, or $\langle{v}|w\rangle$, which is not however measurable; it
is the squared norm
$\left|\langle{v}|w\rangle\right|^2 = \langle{v}|w\rangle \langle{w}|v\rangle$
of the probability amplitude which yields measurable
probabilities.  One of the basic observations which leads to the Jordan
formalism is that these probabilities can be expressed entirely in
terms of the Jordan product of projection operators, since
\begin{equation}
\label{TraceID}
\langle{v}|w\rangle \langle{w}|v\rangle
  = (v^\dagger w)(w^\dagger v)
  \equiv \tr(vv^\dagger \circ ww^\dagger)
\end{equation}
where $\circ$ denotes the \textit{Jordan product}~\cite{Jordan,JNW}
\begin{equation}
\label{JProd}
\AAA \circ \BBB = {1\over2} (\AAA\BBB + \BBB\AAA)
\end{equation}
which is commutative but not associative.

Remarkably, the Jordan formulation of quantum mechanics does not require ($v$
and) $\AAA$ to be complex, but only that the \textit{Jordan identity}
\begin{equation}
\label{JID}
(\AAA\circ \BBB)\circ \AAA^2 = \AAA\circ \left(\BBB\circ \AAA^2\right)
\end{equation}
hold for two Hermitian matrices $\AAA$ and $\BBB$.  As shown
in~\cite{JNW}, the Jordan identity~(\ref{JID}) is equivalent to power
associativity, which ensures that arbitrary powers of Jordan matrices ---
and hence of quantum mechanical observables --- are well-defined.

The Jordan identity~(\ref{JID}) is the defining property of a \textit{Jordan
algebra}~\cite{Jordan}, and is clearly satisfied if the operator algebra is
associative, which will be the case if the elements of the Hermitian matrices
$\AAA$, $\BBB$ themselves lie in an associative algebra.  Remarkably, one
further possibility exists, for which the elements of the Hermitian matrices do
\emph{not} lie in an associative algebra.  This example is the Albert algebra
(also called the exceptional Jordan algebra) $H_3(\OO)$ of $3\times3$
octonionic Hermitian matrices~\cite{JNW,Albert}.
\footnote{The $2\times2$ octonionic Hermitian matrices $H_2(\OO)$ also form a
Jordan algebra, but, even though the octonions are not associative, it is
possible to find an associative algebra that leads to the same Jordan
algebra~\cite{JNW,Jacobson}.}
In what follows we will restrict our attention to this exceptional case.
\footnote{In this case, the equivalence in~(\ref{TraceID}) fails; it is the
right-hand side which provides the correct generalization.}

\section{Quaternions and Octonions}

\begin{figure}
\begin{center}
\includegraphics[height=2in,bb=179 258 433 544]{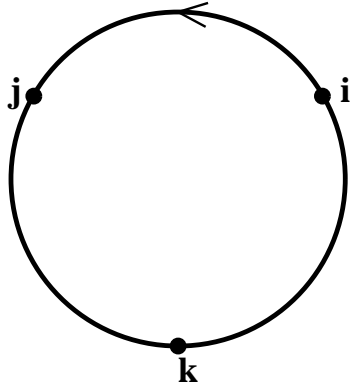}
\end{center}
\caption{The quaternionic multiplication table.}
\label{qmult}
\end{figure}

\begin{figure}
\begin{center}
\includegraphics[height=3in,bb=84 206 560 664]{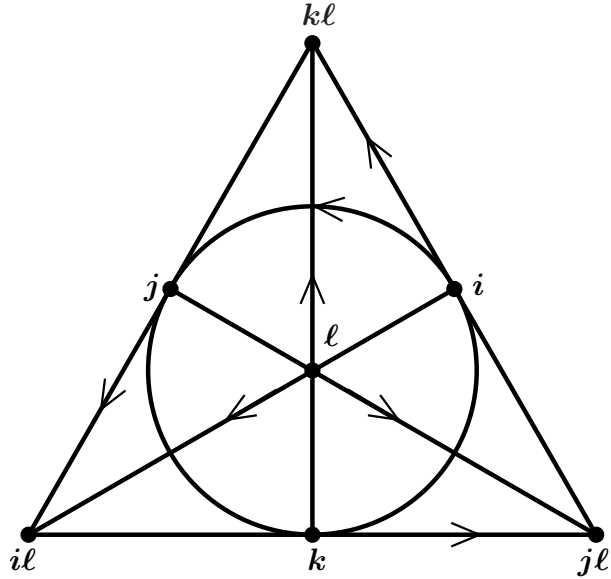}
\end{center}
\caption{The octonionic multiplication table.  Each line segment should be
thought of as circle, identical to the quaternionic multiplication table in
Figure~\ref{qmult}.}
\label{omult}
\end{figure}

The Hurwitz Theorem states that the real numbers $\RR$, complexes $\CC$,
quaternions $\HH$, and octonions $\OO$ are the \textit{only} (normed) division
algebras (over the real numbers).
\footnote{A \textit{division algebra} is a \textit{vector space} over a
\textit{field} (in this case $\RR$) which is also a \textit{ring with
identity} under multiplication, and in which $ax=b$ can be uniquely solved for
$x$ (unless $a=0$).  A \textit{normed} division algebra satisfies~(\ref{norm})
in addition, and is therefore also an \textit{integral domain}, that is, a
ring in which $ab=0$ implies $a=0$ or $b=0$.}
The quaternions and octonions are extensions of the familiar real and complex
numbers.  A quaternion is an arbitrary real linear combination of the real
identity element $1$ and three different square roots of minus one, which are
conventionally called $\{i,j,k\}$ and satisfy the multiplication table given
in Figure~\ref{qmult}.
Similarly, the octonions are formed from seven square roots of minus one which
we will call $\{i,j,k,k\ell,j\ell,i\ell,\ell\}$, whose multiplication table is
summarized in Figure~\ref{omult}.  In these multiplication tables, each point
corresponds to an imaginary unit.  Each line or circle corresponds to a
quaternionic triple with the arrow giving the orientation.  For example,
\begin{eqnarray}
k~\ell &=& k\ell \\
\ell~k\ell &=& k \\
k\ell~k &=& \ell
\end{eqnarray}
and each of these products anticommutes, that is, reversing the order
contributes a minus sign.

We define the \textit{conjugate} $\bar{a}$ of a quaternion or octonion $a$ as
the (real) linear map which reverses the sign of each imaginary unit.  Thus,
if
\begin{equation}
a = a_1+a_2\,i+a_3\,j+a_4\,k+a_5\,k\ell+a_6\,j\ell+a_7\,i\ell+a_8\,\ell
\end{equation}
then
\begin{equation}
\bar{a} = a_1-a_2\,i-a_3\,j-a_4\,k-a_5\,k\ell-a_6\,j\ell-a_7\,i\ell-a_8\,\ell
\end{equation}
Direct computation shows that
\begin{equation}
\bar{ab} = \bar{b}\>\bar{a}
\end{equation}
The \textit{norm} $|a|$ of an octonion $a$ is defined by
\begin{equation}
|a|^2
  = a\bar{a}
  = a_1^2 + a_2^2 + a_3^2 + a_4^2 + a_5^2 + a_6^2 + a_7^2 + a_8^2
\end{equation}
The only octonion with norm $0$ is $0$, and every nonzero octonion
has a unique inverse, namely
\begin{equation}
a^{-1} = {\bar{a} \over ~|a|^2}
\end{equation}
For all the normed division algebras, the norm satisfies the identity
\begin{equation}
|ab| = |a| |b|
\label{norm}
\end{equation}

A remarkable property of the octonions is that they are \textit{not}
associative!  For example, compare
\begin{eqnarray}
(i~j) (\ell) &=& +(k)(\ell) ~=~ +k\ell \\
(i) (j~\ell) &=& (i)(j\ell) ~=~ -k\ell
\end{eqnarray}
However, the octonions are \textit{alternative}, that is, products involving
no more than 2 independent octonions do associate.
The \textit{commutator} of two octonions $a$, $b$ is given as usual by
\begin{equation}
[a,b] = ab - ba
\end{equation}
and we define the \textit{associator} of three octonions $a$, $b$, $c$ by
\begin{equation}
[a,b,c] = (ab)c - a(bc)
\end{equation}
which quantifies the lack of associativity.
More generally, both the commutator and associator are \textit{antisymmetric},
that is, interchanging any two arguments changes the result by a minus sign;
replacing any argument by its conjugate has the same effect, because the real
parts don't contribute to the associator.

The units $i$, $j$, $k$, $k\ell$, $j\ell$, $i\ell$, and $\ell$ are by no means
the only square roots of $-1$.  Rather, \textit{any} pure imaginary quaternion
or octonion squares to a negative number, so it is only necessary to choose
its norm to be $1$ in order to get a square root of $-1$.  The imaginary
quaternions of norm~$1$ form a 2-sphere in the 3-dimensional space of
imaginary quaternions.  The imaginary octonions of norm~$1$ form a 6-sphere in
the 7-dimensional space of imaginary octonions.

Any such unit imaginary quaternion or octonion $\shat$ can be used to
construct a complex subalgebra of $\HH$ or $\OO$, which we will also denote by
$\CC$, and which takes the form
\begin{equation}
\CC = \{a_R + a_s\,\shat\}
\end{equation}
with $a_R,a_s\in\RR$.  Regarding $\shat$ as the complex unit, we have the
familiar Euler identity
\begin{equation}
e^{\shat\theta} = \cos\theta + \shat\sin\theta
\end{equation}
so that \textit{any} quaternion or octonion can be written in the form
\begin{equation}
a = r e^{\shat\theta}
\end{equation}
where
\begin{equation}
r=|a|
\end{equation}
Any \textit{two} unit imaginary octonions $\shat$ and $\that$ that point in
independent directions determine a \textit{quaternionic} subalgebra of $\OO$.

\section{The Structure of \boldmath$G_2$ and \boldmath$SU(3)$}

The freedom to choose an entire \hbox{2-sphere} or \hbox{6-sphere} of square
roots of minus one within the \hbox{3-dimensional} space of the pure imaginary
quaternions or the \hbox{7-dimensional} space of pure imaginary octonions
leads one to investigate the transformations that preserve the corresponding
multiplication table.  These transformations form the automorphism group of
the corresponding division algebra.

In the case of the quaternions, one can imagine rotating $i$ to any pure
imaginary point on the \hbox{2-sphere} (2 degrees of freedom).  Then $j$ can
be chosen to be any direction perpendicular to the direction of $i$, i.e.~on
the equator of the resulting \hbox{2-sphere} (1 degree of freedom).  The
direction of $k$ is determined by the multiplication table.  The
\hbox{3-dimensional} automorphism group of the quaternions is therefore seen
to be $SO(3)$.

For the octonions, one can again imagine rotating $i$ to any pure imaginary
point on the \hbox{6-sphere} (6 degrees of freedom).  Then $j$ must again be
perpendicular to $i$ (5 degrees of freedom) and the direction of $k$ is fixed
by the multiplication table.  But $\ell$ is now free to be any direction
perpendicular to all of the $i$, $j$, and $k$ directions (3 degrees of
freedom) and the directions of the remaining units are determined by the
multiplication table.  This \hbox{14-dimensional} Lie group turns out to be
the exceptional group $G_2$.

\begin{figure}
\begin{center}
\includegraphics[height=3in,bb=78 188 606 664]{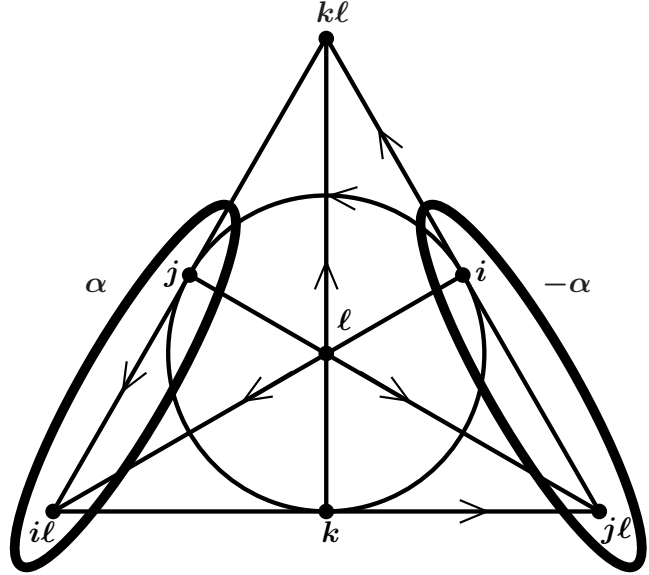}
\end{center}
\caption{One class of elements of $G_2$.  These transformations are also
contained in the preferred $SU(3)$ that fixes $\ell$.}
\label{flip1}
\end{figure}

Another way of envisioning the transformations in the group $G_2$ was first
shown to us by Sudbery~\cite{Tony}.  Consider the octonionic unit $k\ell$ at
the top of the multiplication table shown in Figure~\ref{flip1}.  There are
three pairs of octonionic units that form quaternionic subalgebras with
$k\ell$, i.e.~$\{j,i\ell\}$, $\{j\ell,i\}$, and $\{k,\ell\}$.  We call these
the pairs that ``point to'' $k\ell$.  If the elements of two of these pairs
are rotated into one another oppositely, for instance, if the
$\{j,i\ell\}$-plane is rotated by an angle $\alpha$, and the
$\{j\ell,i\}$-plane is rotated by the angle $-\alpha$, then it turns out that
the multiplication table is preserved.  We have thus constructed a 1-parameter
family of automorphisms.  There are three ways of pairing up the three pairs
of units in this way, but only two are independent.  Since there are 7
different units that can be pointed to, the dimension of this group is again
14.

\begin{figure}
\begin{center}
\includegraphics[height=3in,bb=78 169 606 664]{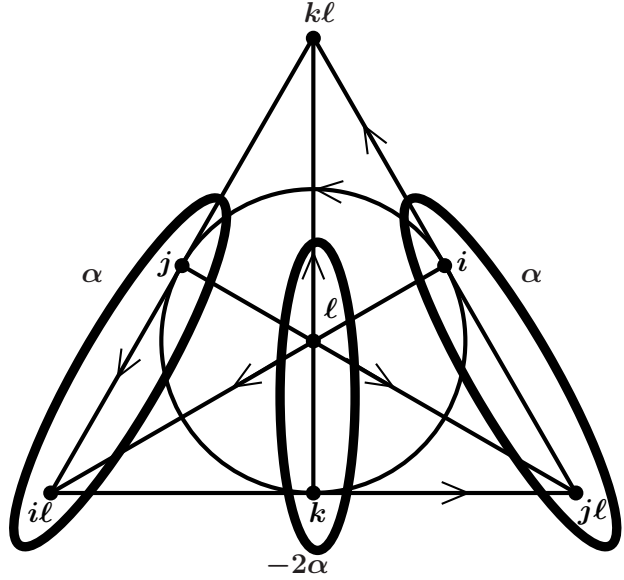}
\end{center}
\caption{A second class of elements of $G_2$.  These transformations are not
contained in the preferred $SU(3)$ that fixes $\ell$.}
\label{flip2}
\end{figure}

\begin{figure}
\begin{center}
\includegraphics[height=3in,bb=78 124 596 793]{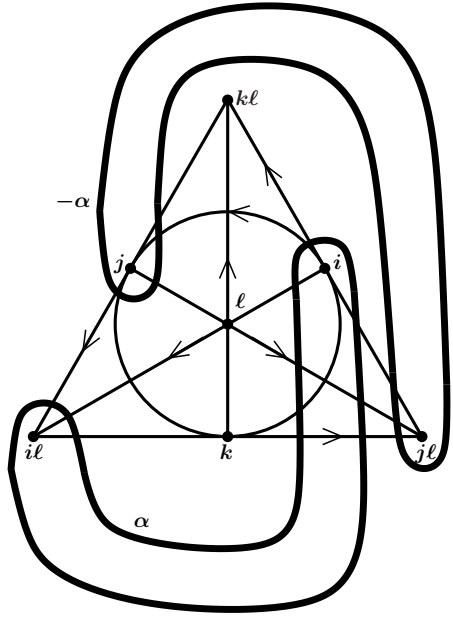}
\end{center}
\caption{A third class of elements of $G_2$.  These transformations are
contained in the preferred $SU(3)$ that fixes $\ell$.}
\label{lfix}
\end{figure}

In what follows, we will need not only $G_2$, but also $SU(3)$, the subgroup
of $G_2$ that fixes one of the octonionic units.  Since $\ell$ is in the
middle of our multiplication table, we will, without loss of generality,
choose it to be the unit that is fixed.  We see that the $G_2$ transformation
in Figure~\ref{flip1} fixes $\ell$ and is therefore in $SU(3)$, but a $G_2$
transformation involving either of the other two pairs that point to $k\ell$
will not fix $\ell$.  To be symmetric, we choose the linear combination of
transformations shown in Figure~\ref{flip2} ($\{j,i\ell\}$ and $\{j\ell,i\}$
both rotate by $\alpha$ and $\{k,\ell\}$ rotates by $-2\alpha$) to be the
$G_2$ transformation that points to $k\ell$ that is \emph{not} in $SU(3)$.  If
we choose to point in turn to each of the six units that are \emph{not}
$\ell$, we have six $G_2$ transformations that are in $SU(3)$ and six that are
not.  What about the remaining two $G_2$ transformations?  These are
transformations that point to $\ell$.  One such transformation, shown in
Figure~\ref{lfix}, rotates $\{i\ell, i\}$ by $\alpha$ and $\{j\ell,j\}$ by
$-\alpha$.  There are three transformations of this type, all of which fix
$\ell$ and are therefore elements of $SU(3)$, but only two are linearly
independent.  Any two of these transformations complete the 8-dimensional Lie
group $SU(3)$.

Yet another way to describe $G_2$ is in terms of \textit{inner automorphisms},
that is, transformations of the form
\begin{equation}
x\mapsto axa^{-1}
\label{conjugate}
\end{equation}
Inner automorphisms always preserve an associative multiplication rule, since
\begin{equation}
(axa^{-1})(aya^{-1}) = a(xy)a^{-1}
\label{inner}
\end{equation}
However, this condition is nontrivial over the octonions, since the
parentheses cannot be moved.  As shown in~\cite{Lorentz},~(\ref{inner}) holds
for all $x,y\in\OO$ if and only if $a$ is a sixth root of unity.  That is, the
inner automorphisms of the octonions are given by~(\ref{conjugate}) where
\begin{equation}
a = e^{n\pi\shat/3}
\label{root}
\end{equation}
where $n\in\ZZ$ and $\shat$ is any pure imaginary unit octonion.  As further
discussed in~\cite{Lorentz}, any $G_2$ transformation can in fact be generated
by a finite sequence of nested transformations of the form~(\ref{conjugate}),
with $a$ given by~(\ref{root}).

\section{The Jordan Eigenvalue Problem}

In previous work~\cite{Other}, we solved the \textit{Jordan eigenvalue
problem}, namely the eigenmatrix problem
\begin{equation}
\label{JEigen}
\AAA \circ \VVV = \lambda \VVV
\end{equation}
where $\AAA$ and $\VVV$ are both $3\times3$ octonionic Hermitian matrices.
Unlike the right eigenvalue problem $\AAA v=v\lambda$ considered
in~\cite{Eigen}, the Jordan eigenvalue problem~(\ref{JEigen}) admits only real
eigenvalues, which do solve the characteristic equation for $\AAA$,
namely~\cite{Freudenthal}
\begin{equation}
\label{Char}
-\det(\AAA-\lambda\,\III)
  = \lambda^3 - (\tr\,\AAA)\,\lambda^2
        + \sigma(\AAA)\,\lambda - (\det \AAA)\,\III
  = 0
\end{equation}
where $\III$ denotes the identity matrix, $\sigma(\AAA)$ is defined by
\begin{equation}
\label{Sigma}
\sigma(\AAA)
  = {1\over2} \left( (\tr\,\AAA)^2 - \tr (\AAA^2) \right)
  = \tr(\AAA*\AAA)
\end{equation}
the operation $*$ denotes the Freudenthal product
\begin{equation}
\AAA*\BBB
  = \AAA \circ \BBB - {1\over2} \Big(\AAA\,\tr(\BBB)+\BBB\,\tr(\AAA)\Big)
        + {1\over2} \Big(\tr(\AAA)\,\tr(\BBB)-\tr(\AAA\circ \BBB)\Big) \,\III
\label{FProd}
\end{equation}
and the determinant is defined unambiguously by
\begin{equation}
\label{Det}
\det(\AAA) = {1\over3} \, \tr \Big( (\AAA*\AAA) \circ \AAA \Big)
\end{equation}
The Jordan and Freudenthal products are generalizations of the standard dot
and cross products.

Just as in the more familiar complex case, normalized eigenmatrices for each
nondegenerate eigenvalue are primitive idempotents, and the degenerate case
can be handled using Gram-Schmidt orthogonalization.  Furthermore, each
primitive idempotent is in fact an element of the \textit{Cayley-Moufang
plane} $\OP^2$, which can be characterized as
\begin{eqnarray}
\OP^2
  &:=& \{\VVV\in H_3(\OO): \VVV\circ\VVV=\VVV \hbox{ and } \tr\VVV=1\}
	\label{op2}\\
  &\equiv& \{\VVV\in H_3(\OO): \VVV*\VVV=0 \hbox{ and } \tr\VVV=1\}\nonumber
\end{eqnarray}
It is straightforward to show from the first condition in (\ref{op2}) that
the components of any element $\VVV\in\OP^2$ must lie in some quaternionic
subalgebra of $\OO$, which of course depends on $\VVV$.  Put differently, the
associator of the (independent, off-diagonal) components of $\VVV$, denoted
$[\VVV]$, must vanish.  But quaternionic primitive idempotents have
(nonunique) ``square roots'', $\VVV=\Psi\Psi^\dagger$, so that we can also
write
\begin{equation}
\OP^2 = \{\Psi\Psi^\dagger:\Psi\in\OO^3, [\Psi]=0, \Psi^\dagger\Psi=1\}
\end{equation}
where $[\Psi]$ denotes the associator of the components of $\Psi$.  We refer
to such 3-component octonionic column vectors $\Psi$ as \textit{Cayley
spinors}.

Putting this all together, any $3\times3$ octonionic Hermitian matrix $\AAA$
can be expressed as the sum of the squares of quaternionic columns, which are
orthogonal under the Jordan product, that is
\begin{equation}
\label{Decomp}
\AAA = \sum_{i=1}^3 \lambda_i \VVV_i
\end{equation}
in terms of primitive idempotents $\VVV_i=\Psi_i\Psi_i^\dagger\in\OP^2$ and
their corresponding eigenvalues $\lambda_i$.

\section{The Structure of \boldmath$E_6$}

The automorphism group of the Jordan product~(\ref{JProd}) (and consequently
also of the Freudenthal product~(\ref{FProd})) is the exceptional group $F_4$,
and the group which leaves the determinant~(\ref{Det}) invariant is a
particular real form of the exceptional group $E_6$.  These groups can
be interpreted as $F_4=SU(3,\OO)$ and $E_6=SL(3,\OO)$, as we now show; for
further details see~\cite{Denver}.

In previous work~\cite{Lorentz}, Manogue and Schray showed how to write
the (double cover of the) Lorentz group $SO(9,1)$ as $SL(2,\OO)$, with the
action given by
\begin{equation}
\XX \longmapsto \MM \XX \MM^\dagger
\end{equation}
where $\XX\in H_2(\OO)$, the $2\times2$ Hermitian matrices with octonionic
components.  The key to that work was to give an explicit set of basis
transformations --- the rotations and boosts in coordinate planes --- which
were \textit{compatible} with the spinor representation in the sense that if
$\theta\in\OO^2$ transforms like
\begin{equation}
\theta \longmapsto \MM \theta
\end{equation}
then there are no associativity problems in the vector transformation
\begin{equation}
\MM (\theta\theta^\dagger) \MM^\dagger = (\MM \theta)(\MM \theta)^\dagger
\end{equation}

Any such basis transformation $\MM\in SL(2,\OO)$ can be immediately
reinterpreted as a $3\times3$ transformation $\MMM$ via
\begin{equation}
\MMM = \pmatrix{\MM& 0\cr 0& 1\cr}
\label{typeI}
\end{equation}
and it is straightforward to verify that any such $\MMM$ preserves the
determinant of $\XXX\in H_3(\OO)$ and is therefore an element of $E_6$.

How many such transformations are there?  We first give the basis
transformations for the simpler case of $SL(2,\CC)$, adapted
from~\cite{Lorentz} and rewritten as elements of $E_6$.  When interpreting
these transformations, it is helpful to recall that, in this case,
\begin{equation}
\XX = \pmatrix{t+z& x-\ell y\cr \noalign{\smallskip}  x+\ell y& t-z\cr}
\label{Vector}
\end{equation}
We have the three rotations
\begin{equation}
\MMM_{xy}
  = \pmatrix{e^{-\ell\theta/2}& 0& 0\cr 0& e^{\ell\theta/2}& 0\cr 0& 0& 1\cr}
\label{Rz}
\end{equation}
\begin{equation}
\MMM_{yz} = \pmatrix{~~\cos\frac\theta2& -\ell\sin\frac\theta2& 0\cr
		\noalign{\smallskip}
                -\ell\sin\frac\theta2& ~~\cos\frac\theta2& 0\cr
		\noalign{\smallskip}
		0& 0& 1\cr}
\qquad
\MMM_{zx} = \pmatrix{\cos\frac\theta2& -\sin\frac\theta2& 0\cr
		\noalign{\smallskip}
                \sin\frac\theta2& ~~\cos\frac\theta2& 0\cr
		\noalign{\smallskip}
		0& 0& 1\cr}
\end{equation}
and the three boosts
\begin{equation}
\MMM_{tz} = \pmatrix{e^{\beta/2}& 0& 0\cr 0& e^{-\beta/2}& 0\cr 0& 0& 1\cr}
\label{tzboost}
\end{equation}
\begin{equation}
\MMM_{tx} = \pmatrix{\cosh\frac\beta2& \sinh\frac\beta2& 0\cr
		\noalign{\smallskip}
                \sinh\frac\beta2& \cosh\frac\beta2& 0\cr
		\noalign{\smallskip}
		0& 0& 1\cr}
\qquad
\MMM_{ty} = \pmatrix{\cosh\frac\beta2& -\ell\sinh\frac\beta2& 0\cr
		\noalign{\smallskip}
                \ell\sinh\frac\beta2& ~~\cosh\frac\beta2& 0\cr
		\noalign{\smallskip}
		0& 0& 1\cr}
\end{equation}
written in terms of the single imaginary unit $\ell$.  To generate a
$3\times3$ representation of $SL(2,\OO)$, start by replacing $\ell$ in turn by
each of the other imaginary units, yielding 18 new transformations, for a
total of 24.  The remaining 21 transformations in $SL(2,\OO)$ are precisely
the rotations of the imaginary units with each other, which make up an $SO(7)$
(really $\hbox{\textit{Spin}}(7)$; we are being casual about double covers).
As shown in~\cite{Lorentz}, these rotations are obtained by \textit{nesting},
that is by transformations of the form
\begin{equation}
\XXX \longmapsto \MMM_2 \left( \MMM_1 \XXX \MMM_1^\dagger \right) \MMM_2^\dagger
\label{nest}
\end{equation}
where each corresponding $\MM$ represents a ``flip'', that is, a pure
imaginary multiple of the ($2\times2$!)\ identity matrix.  Thus, a typical
$\MMM$ takes the form
\begin{equation}
\MMM =  \pmatrix{\ell& 0& 0\cr 0& \ell& 0\cr 0& 0& 1\cr}
\end{equation}
where it is important to note that $\MMM$ is \emph{not} a multiple of the
($3\times3$) identity matrix.

We are now ready to count the basis transformations of $E_6$.  At first sight,
it appears we have three copies of $SL(2,\OO)$ --- simply repeat the
embedding~(\ref{typeI}) with the two other obvious block structures.  We call
these three copies type I, II, and III.  However, this yields $3\times45=135$
transformations, and, while these transformations do indeed generate all of
$E_6$, it is clear that they can not be a basis, since the dimension of $E_6$
is only $78$.

Let's try again.  Each of these three copies of $SL(2,\OO)=SO(9,1)$ contains a
copy of $SO(8)$.  A famous property of $SO(8)$ called \textit{triality}
asserts in this context that these three copies of $SO(8)$ in fact consist of
the same $E_6$ transformations (but labeled differently), so we should count
these copies only once.  But $SO(8)$ has 28 elements, to which we must add 3
copies of the $8$ rotations needed to get to $SO(9)$, then 3 copies of the $9$
boosts needed to get to $SO(9,1)$, resulting in $28+3\times8+3\times9=79$
transformations.  A bit of thought reveals that the 3 copies of the
$tz$-boost~(\ref{tzboost}) are not independent; removing one of them correctly
yields an explicit set of $78$ basis transformations for $E_6$, also
justifying the interpretation $E_6=SL(3,\OO)$.

It is worth pointing out that, due to triality, only the 14 $G_2$
transformations need to be written in the nested form~(\ref{nest}).
Remarkably, the remaining 14 $SO(8)$ transformations can all be expressed
using the type I transformation~(\ref{Rz}) and its types II and III variants.
The former are just the usual 7 rotations needed to get from $SO(7)$ to
$SO(8)$, but the latter yield an unnested description of the 7 non-$G_2$
transformations in $SO(7)$, which take the form
\begin{equation}
\MMM_\ell = \pmatrix{e^{\ell\theta/2}& 0& 0\cr
                        0& e^{\ell\theta/2}& 0\cr 0& 0& e^{-\ell\theta}\cr}
\label{Rl}
\end{equation}
Each of these transformations rotates 3 octonionic planes by the same amount,
and can therefore be thought of as a ``phase'' transformation.

What about $F_4$?  Note that we have described $27-1=26$ boosts, and
$78-26=52$ rotations.  So our $E_6$ is the real representation with 26 boosts,
commonly written as $E_{6(-26)}$, with the number in parentheses denoting the
number of boosts minus the number of rotations.  It is straightforward to show
that $F_4$ preserves the trace of elements of $H_3(\OO)$, corresponding to the
timelike direction; this is the compact representation of $F_4$, consisting
precisely of the rotation subgroup of this real form of $E_6$.  These
considerations justify the interpretation $F_4=SU(3,\OO)$.

Returning to the characteristic equation~(\ref{Char}), not only does $E_6$
preserve the determinant, it also preserves the condition $\sigma=0$.  But
these two coefficients control the number of nonzero eigenvalues --- 3 if
$\det\AAA\ne0$, 2 if $\det\AAA=0\ne\sigma(\AAA)$, and 1 if
$\det\AAA=0=\sigma(\AAA)$.  Thus, $E_6$ preserves the number of nonzero
eigenvalues of $\AAA$, and hence the number of terms in the
decomposition~(\ref{Decomp}) with nonzero eigenvalue.

\section{Symmetry Breaking and Particle Physics}
\label{Breaking}

In order to apply this formalism to elementary particle physics, we break the
full symmetry and explore how a Jordan matrix transforms under various
subgroups of $E_6$.  One such symmetry breaking occurs when we choose a
preferred $SL(2,\OO)$ subgroup of $SL(3,\OO)$, as in~(\ref{typeI}).  This
leads us to impose a block structure on $H_3(\OO)$, so that
\begin{equation}
\PPP = \pmatrix{\PP& \psi\cr \psi^\dagger& n\cr}
\end{equation}
where $\PP\in H_2(\OO)$ transforms like a 10-dimensional momentum vector,
$\psi\in\OO^2$ transforms like a (Majorana-Weyl) spinor, and $n\in\RR$ is a
scalar.  Direct computation shows that
\begin{equation}
\PPP*\PPP
  = \pmatrix{\tilde{\psi\psi^\dagger}-n\tilde{\PP}& \tilde{\PP}\psi\cr
                \noalign{\smallskip}
                (\tilde{\PP}\psi)^\dagger& \det \PP\cr}
\label{Star}
\end{equation}
where tilde denotes trace reversal, that is, $\tilde\PP=\PP-\tr(\PP)\,\II$.
As shown in~\cite{Spin}, the massless, momentum-space Dirac equation in 10
dimensions can be written
\begin{equation}
\tilde{\PP} \psi = 0
\label{Dirac}
\end{equation}
which implies the nonlinear constraint
\begin{equation}
\det \PP = 0
\label{Massless}
\end{equation}
The general solution of (\ref{Dirac}) and (\ref{Massless}) is of the form
\begin{eqnarray}
\psi &=& \theta\xi \\
\PP &=& \theta\theta^\dagger
\end{eqnarray}
where the components of $\theta\in\OO^2$ lie in the complex subalgebra of
$\OO$ determined by $\PP$ and $\xi\in\OO$ is arbitrary.  Using~(\ref{Star}),
these equations are seen to be precisely the same as
\begin{equation}
\PPP*\PPP = 0
\label{DiracStar}
\end{equation}
where
\begin{equation}
n = |\xi|^2
\end{equation}
Thus, (normalized) solutions of the Dirac equation are precisely elements of
$\OP^2$, and therefore the squares of Cayley spinors.

In previous work~\cite{Spin}, we discussed solutions of the Dirac equation in
the form~(\ref{Dirac}).  Remembering that solutions of~(\ref{DiracStar}) are
quaternionic, and reducing 10 spacetime dimensions to 4 by the simple
expedient of choosing a preferred complex subalgebra of $\OO$, we used spin
eigenstates and particle/antiparticle projection operators as usual to
identify particle states within our division algebra formalism.  This
procedure identified a spin-$\frac12$ massive particle with two spin states,
namely
\begin{eqnarray}
\eplus\phd &=& \pmatrix{1\cr k\cr} \qquad
  \eplus\phd \eplus^\dagger = \pmatrix{1&-k\cr k&~~1\cr}\\
\eminus\phd &=& \pmatrix{-k\cr ~~1\cr} \qquad
  \eminus\phd \eminus^\dagger = \pmatrix{1&-k\cr k&~~1\cr}
\end{eqnarray}
where the direction of the arrow indicates the $z$-component of the spin, and
where the second equality in each case gives the momentum vector.  Comparison
with~(\ref{Vector}) shows that the $x$, $y$, and $z$ components of the
momentum vanish; these states are given at rest.  Similarly, there is an
analogous antiparticle with two spin states.  The procedure also identified a
left-handed massless particle, which when moving in the $z$-direction takes
the form
\begin{eqnarray}
\nu_z\phd &=& \pmatrix{0\cr k\cr} \qquad
  \nu_z\phd \nu_z^\dagger = \pmatrix{0&0\cr 0&1\cr}
\end{eqnarray}
together with a massless particle of the opposite helicity, which when moving
in the $z$-direction takes the form
\begin{equation}
\hbox{\O}_z\phd = \pmatrix{0\cr 1\cr} \qquad
  \hbox{\O}_z\phd \hbox{\O}_z^\dagger = \pmatrix{0&0\cr 0&1\cr}
\end{equation}
Any other 2-component quaternionic column can be identified as an
appropriately rotated and/or boosted superposition of these particles in the
usual way.  Noting that the first two particles carry an octonionic label
($k$), it is straightforward to generalize these particles to 3 generations of
leptons, labeled by $i$, $j$, $k$; the remaining particle does not have an
octonionic label.  Thus, an octonionic description of the Dirac equation in 10
dimensions yields a particle spectrum containing precisely 3 generations of
leptons, each with a single-helicity, massless neutrino, together with a
single ``sterile'' neutrino of the opposite helicity but with no generation
structure.

\section{Discussion}

We have shown how to break the symmetry group $E_6$ so that a Lorentzian $3+1$
dimensional momentum space emerges, together with internal symmetries that
describe the correct spin/helicity transformations on Cayley spinors to
describe leptons.  Furthermore, precisely three generations of such leptons
exist which respect the octonionic structure of the transformations.
Contained naturally within this description of leptons and their symmetries
are three massless left-handed neutrinos and a single, sterile, right-handed
neutrino.

\begin{figure}
\begin{center}
\includegraphics[bb=84 353 478 685]{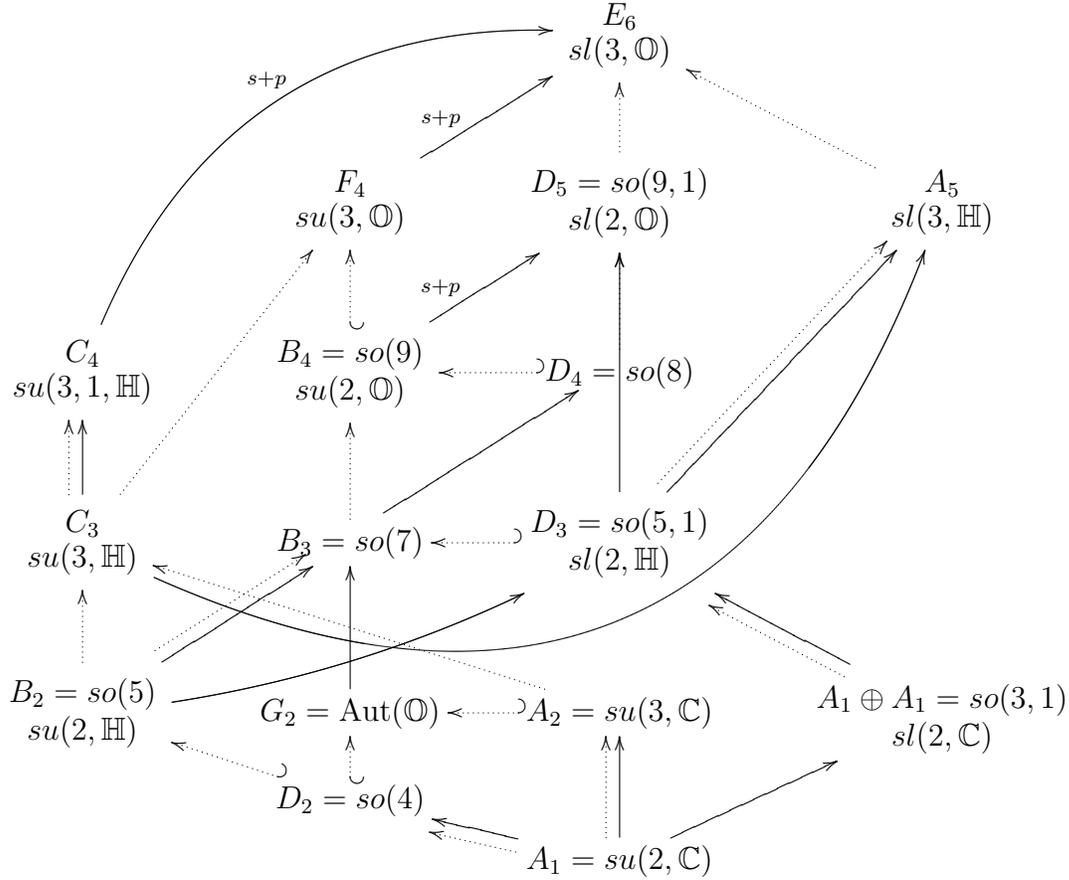}
\end{center}
\caption{A map of $E_6$, taken from~\cite{Aaron}.}
\label{E6map}
\end{figure}

We have learned several important lessons along the way.  First, in spite of
the non-commutativity and non-associativity of the octonions, everything that
one might want to do can be made to work if one defines everything carefully.
Second, when working with the octonions, it is important to make the Lie group
structure primary, rather than the Lie algebra structure.  Some of the nested
group transformations described in (\ref{nest}) cannot be described in terms
of the exponentials of any Lie algebra transformations.  Finally, in order to
keep a complex structure on the Lie algebra from interfering with the
octonionic units in the matrices, it is important in the symmetry-breaking
process to look at \emph{real} Lie subalgebras rather than complexified ones.
We call the reader's attention to the recent work of Aaron
Wangberg~\cite{Aaron}, in which the important real forms of the subgroups of
$E_6$ were identified, using a division-algebra perspective.  A map of these
subgroups, taken from~\cite{Aaron}, appears in Figure~\ref{E6map}.

\goodbreak
We conclude with several speculative questions:
\begin{enumerate}
\item
In this paper, we have described leptons in terms of ``1-squares,'' i.e.~the
squares of Cayley spinors.  But the most general octonionic hermitian matrix
can be a linear combination of up to three such squares.  Is it possible to
describe the meson and baryon sectors of standard particle physics as
``2-squares'' and ``3-squares,'' respectively?
\item
Does the group $E_6$, broken in the way we have discussed, describe both the
standard model interactions and Lorentz transformations in $3+1$-dimensions?
\item
After the symmetry breaking described in Section~\ref{Breaking}, some of the
Lie group transformations of the form~(\ref{conjugate}),~(\ref{root}) are no
longer connected to the identity.  Do these discrete transformations
correspond to discrete conserved quantities such as charge?
\end{enumerate}

\ack
We would like to acknowledge the many contributions of colleagues over the
years who have contributed to our understanding of this field: David Fairlie,
Ed Corrigan, J\"org Schray, Robin Tucker, Jason Janesky, Susumu Okubo, Aaron
Wangberg, Jim Wheeler, and, of course, Tony Sudbery.  CAM thanks the
organizers for support to attend the TonyFest and present this work.  We thank
the Department of Physics at Utah State University for support and hospitality
during the preparation of this manuscript.  This work was supported in part by
a grant from the Foundational Questions Institute (FQXi).

\goodbreak


\section*{References}
\begin{thebibliography}{99}

\bibitem{Sudbery}
Sudbery A and Manogue C A 1989
General solutions of the covariant superstring equations of motion
\textit{Phys. Rev. D} \textbf{40} 4073--4077

\bibitem{Jordan}
Jordan P 1933
\"Uber die Multiplikation quantenmechanischer Gr\"o\ss en
\textit{Z. Phys.} \textbf{80} 285--291

\bibitem{JNW}
Jordan P, von Neumann J and Wigner E 1934
On an algebraic generalization of the quantum mechanical formalism
\textit{Ann. of Math.} \textbf{35} 29--64

\bibitem{GT}
G\"ursey F and Tze C-H 1996
\textit{
  On the Role of Division, Jordan, and Related Algebras in Particle Physics}
(Singapore: World Scientific)

\bibitem{Other}
Dray T and Manogue C A 1999
The Exceptional Jordan Eigenvalue Problem
\textit{Internat. J. Theoret. Phys.\/} \textbf{38} 2901--2916
(\textit{Preprint} \url{math-ph/9910004})

\bibitem{Albert}
Albert A A 1934
On a certain algebra of quantum mechanics
\textit{Ann. Math.} \textbf{35} 65--73

\bibitem{Jacobson}
Jacobson N 1968
\textit{Structure and Representations of Jordan Algebras}
(\textit{Amer.  Math.  Soc. Colloq. Publ.} \textbf{39})
(Providence: American Mathematical Society)

\bibitem{Tony}
Sudbery A 1988
private communication

\bibitem{Lorentz}
Manogue C A and Schray J 1993
Finite Lorentz transformations, automorphisms, and division algebras
\textit{J. Math. Phys.} \textbf{34}, 3746--3767
(\textit{Preprint} \url{hep-th/9302044})

\bibitem{Eigen}
Dray T and Manogue C A 1998
The octonionic eigenvalue problem
\textit{Adv. Appl. Clifford Algebras} \textbf{8}, 341--364
(\textit{Preprint} \url{math/9807126})

\bibitem{Freudenthal}
Freudenthal H 1953
Zur ebenen Oktavengeometrie
\textit{Proc. Kon. Ned. Akad. Wet.} \textbf{A56}, 195--200

\bibitem{Denver}
Dray T and Manogue C A 2010
Octonionic Cayley Spinors and $E_6$
\textit{Comment. Math. Univ. Carolin.} to appear
(\textit{Preprint} \url{arXiv:0911.2255})

\bibitem{Aaron}
Wangberg A 2007
\textit{The Structure of $E_6$} Ph.D. thesis Oregon State University
(\textit{Preprint} \url{arXiv:0711.3447})

\bibitem{Spin}
Manogue C A and Dray T 2000
Quaternionic spin
\textit{Clifford Algebras and their Applications in Mathematical Physics}
eds R Ab\l amowicz and B Fauser
(Boston: Birkh\"auser) pp 29--46
(\textit{Preprint} \url{ hep-th/9910010})

\end{thebibliography}

\end{document}